\documentclass{article}
\usepackage[utf8]{inputenc}
\usepackage{amsmath, amsthm, amssymb, amsfonts}
\usepackage{bm}
\usepackage{dsfont}
\usepackage[utf8]{inputenc}
\usepackage{rotate}
\usepackage{tikz}
\usepackage{tikz-cd}
\usepackage[arrow,matrix,curve]{xy} 	
\usepackage{xcolor}

\theoremstyle{plain}
\newtheorem{theorem}{Theorem}

\newtheorem{corollary}[theorem]{Corollary}

\theoremstyle{definition}

\newtheorem{remark}[theorem]{Remark}
\newtheorem{example}[theorem]{Example}

\newcommand{\norm}[1]{\left\lVert#1\right\rVert}

\title{Approximation of Discrete Measures by Finite Point Sets}
\author{Christian Wei\ss{}}
\date{\today}

\begin{document}

\maketitle

\begin{abstract}
    For a probability measure $\mu$ on $[0,1]$ without discrete component, the best possible order of approximation by a finite point set in terms of the star-discrepancy is $\frac{1}{2N}$ as has been proven relatively recently. However, if $\mu$ contains a discrete component no non-trivial lower bound holds in general because it is straightforward to construct examples without any approximation error in this case. This might explain, why the approximation of discrete measures on $[0,1]$ by finite point sets has so far not been completely covered in the existing literature. In this note, we close this gap by giving a complete description of the discrete case. Most importantly, we prove that for any discrete measure the best possible order of approximation is for infinitely many $N$ bounded from below by $\frac{1}{cN}$ for some constant $c \geq 2$ which depends on the measure. This implies, that for a finitely supported discrete measure on $[0,1]^d$ the known possible order of approximation $\frac{1}{N}$ is indeed the optimal one.
\end{abstract}

\section{Introduction}
According to \cite{FGW21}, the Lebesgue measure is the hardest Borel measure on $[0,1]$ to approximate by a finite point set. In order to formulate the result in a mathematically precise way, recall first that the star-discrepancy between two probability measures $\mu, \nu$ on $[0,1]$ is defined by
$$
D^*_{N}(\mu;\nu) := \sup_{A \in \mathcal{A}} \left| \mu(A) - \nu(A)\right|,
$$
where $\mathcal{A}$ is the set of all half-open intervals in $[0,1]$ which have one vertex at the origin. Furthermore, the probability measure associated to a finite set $(y_i)_{i=1}^N$ is given by
\begin{equation} \label{e:nun}
	\nu_{N} = \frac{1}{N} \sum_{i=1}^{N} \delta_{y_{i}},
\end{equation}
where $\delta_{y_{i}}$ denotes the Dirac measure centered at $y_i$. Borel measures on the interval $[0,1]$ have a particularly comprehensible structure. Lebesgue's decomposition theorem states that any Borel measure $\mu$ can be written as
	$$\mu = \mu_{ac} + \mu_{d} + \mu_{cs},$$
where $\mu_{ac}$ is absolutely continuous with respect to the Lebesgue measure, that is $\mu_{ac}$ is zero on sets of Lebesgue measure zero, $\mu_d$ is a discrete measure, that is, it is zero on the complement of some countable set, and $\mu_{cs}$ is continuous singular, that is, $\mu_{cs}$ is zero on the complement of some set $B$ of Lesbesgue measure zero but assigns no weight to any countable set of points, see e.g. \cite{HS75}, Chapter V. Based on this observation, the following result holds.
\begin{theorem}[\cite{FGW21}, Theorem 1.5] \label{thm:01}
	Fix $\mu$ a probability measure on $[0,1]$. 
	\begin{itemize}
		\item[(i)] For all $N\in \mathbb{N}$, there exists a finite set $(y_i)_{i=1}^N$ such that $\nu_N$ as in \eqref{e:nun} satisfies 
		$$D^*_{N}\left(\mu; \nu_N \right) \leq \frac{1}{N}.$$
		\item[(ii)] Suppose $\mu$ is a probability measure with no point masses. That is, $\mu = \mu_{ac} + \mu_{cs}$. Then 
		\begin{equation}
			D_N^*\left(\mu ; \nu_N \right) \geq \frac{1}{2N}
		\end{equation}
		for any finite set $(y_i)_{i=1}^N$ and $\nu_{N}$ as in \eqref{e:nun}. 
	\end{itemize}
\end{theorem}
This result answered the general question from \cite{ABN18}, where the authors asked for arbitrary dimensions which Borel measure on $[0,1]^d$ is the hardest to approximate by finite point sets,  in the one-dimensional case. However, the lower bound in Theorem~\ref{thm:01} (ii), i.e. the best theoretically possible speed of approximation by finite point sets, is restricted to measures without discrete components only. The main purpose of this note is to close this gap and thereby to complete the discussion on approximation of measures by finite point sets in the one-dimensional case.\\[12pt]
Already the simplest possible example of approximating the Dirac measure centered at $x_0 \in [0,1]$ by a finite point set yields some insight: indeed, it is possible in the discrete case (in contrast to the other cases) to have an approximation error $D_N^*(\mu,\nu_N)=0$. However, this is a rather special situation as our main theorem shows.
\begin{theorem} \label{main:thm} Fix a discrete probability measure $\mu$ on $[0,1]$.
\begin{itemize}
    \item[(i)] \textbf{Finitely supported, rational weights:} Let $\mu$ be given by
    $$\sum_{i=1}^n \xi_i x_i$$
    for $0 \leq x_1 < x_2 < \ldots < x_n \leq 1$, where $\xi_i = \frac{p_i}{q_i}$ with $p_i,q_i \in \mathbb{N}$ and $\gcd(p_i,q_i) = 1$. If $\textrm{lcm}(q_1,\ldots,q_{n-1})|N$, then there exists a finite point set $(y_k)_{k=1}^N$ with $D_N^*(\mu,\nu_N) = 0$. Otherwise $D_N^*(\mu,\nu_N) \geq \frac{1}{qN}$ holds for all finite point sets $(y_k)_{k=1}^N$, where $q = \max(q_1,\ldots,q_n)$. 
    \item[(ii)] \textbf{Finitely supported, irrational weights:} Let $\mu$ be given by
    $$\sum_{i=1}^n \xi_i x_i$$
    for $0 \leq x_1 < x_2 < \ldots < x_n \leq 1$ and $\xi_i \notin \mathbb{Q}$ with $\sum_{i=1}^l \xi_i \notin \mathbb{Q}$ for all $l=1,\ldots,n-1$. Then for any $c>2$, there exist infinitely many $N$ such that
    $$D_N^*(\mu,\nu_N) \geq \frac{1}{cN}$$
    holds for any finite point set $(y_k)_{k=1}^N$.
    \item[(iii)] \textbf{Infinitely supported measures:} Let $\mu$ be given by
    $$\sum_{i=1}^\infty \xi_i x_i$$
    with $x_i < x_{i+1}$ for all $i \in \mathbb{N}$ and $\xi_i > 0$ for all $i \in \mathbb{N}$. Then there exists a constant $c \geq 2$ such that
    $$D_N^*(\mu,\nu_N) \geq \frac{1}{cN}$$
    holds for infinitely many $N$ and any finite point set $(y_k)_{k=1}^N$. If in addition $\sum_{1=1}^l \xi_i \notin \mathbb{R}$ for all $l \in \mathbb{N}$, then $c > 2$ can be chosen arbitrarily.
\end{itemize}
\end{theorem}
\begin{remark} If in (ii) the condition $\sum_{i=1}^l \xi_i \notin \mathbb{Q}$ for all $l=1,\ldots,n$ is violated, then the situation can be treated similarly as in case (i) and we would need to consider the denominators of $\sum_{i=1}^l \xi_i$ to derive a lower bound for infinitely many $N$.
\end{remark}
It is possible to use our approach also in higher dimensions and we again obtain a lower bound of the form $D_N^*(\mu,\nu_N) \geq \frac{1}{cN}$ for infinitely many $N$. Together with \cite{FGW21}, Proposition~2.2, this leads to the following interesting corollary.
\begin{corollary} Let $\mu$ be a probability measure on $[0,1]^d$ which is supported on a finite number of points $k \in \mathbb{N}$. Then there exists constant $c_{\mu}$, which depends on the measure, such that
$$D_N^*(\mu,\nu_N) \geq \frac{1}{c_{\mu}N}$$
for infinitely $N$ and arbitrary point sets $(y_i)_{i=1}^N$. Moreover for any $N \in \mathbb{N}$, there exists a constant $C_k$, which only depends on the number of points, and a finite set $(y_i)_{i=1}^N$ such that
$$D_N^*(\mu,\nu_N) \leq \frac{C_k}{N}.$$
\end{corollary}
In other words, in any dimension any finitely supported discrete measure can be approximated by a finite set of order $\frac{1}{N}$ and this is the best possible order of approximation. If the decay rate of the weights is fast enough, it is furthermore known due to \cite{FGW21}, Theorem 1.1, that certain measures can be approximated by order of convergence at most $\log(N)/N$. Nonetheless, we expect that due to the combinatorial richness of inclusion of half-open intervals in higher dimensions, there exist infinitely supported discrete measures in dimensions $d \geq 3$ with a bigger minimal possible order of approximation. We therefore ask the question under which conditions on the probability measure the lower bounds from Theorem~\ref{main:thm} are also optimal in higher dimensions. This question is of particular interest because lower bounds for the star-discrepancy, e.g. for the Lebesgue measure, are typically very hard obtained, compare \cite{KN74,Nie92}.\\[12pt]
The main reasons why Theorem~\ref{main:thm} holds can be easiest understood by considering the second simplest example, namely discrete measures supported on two points $x_1,x_2$ only. It turns out that \textbf{Kronecker sequences}, which are for $\alpha \in \mathbb{R}$ defined by $\{i \alpha\}_{i=1}^\infty$, where $\{ \cdot \}$ denotes the fractional part of a real number (see e.g. \cite{KN74}), appear here prominently. 
\begin{example} \label{ex:n=2} As an illuminating example let us start with the case that $\mu$ consists of two point masses only, $n=2$. For fixed $N \in \mathbb{N}$, it is clear that some of the weight, i.e. some of the $y_i$, needs to be placed at $x_1$ and the rest at $x_2$. At first, we consider the case $\xi_1,\xi_2 \in \mathbb{Q}$ and let $\xi_i = p_i/q_i$ with $\gcd(p_i,q_i)=1$ for $i=1,2$. If $q_1|N$, then also $q_2|N$. Choosing $k_1 = p_1 \cdot N /q_1$ times $y_i = x_1$ and $k_2 = p_2 \cdot N / q_2$ times $y_i = x_2$ yields an approximation error $D_N^*(\mu,\nu_N)=0$ as predicted by Theorem~\ref{main:thm}. If on the other hand $q_1 \nmid N$, write  $q_1 = r_1 s_1$ with $r_1 | N$ and $\gcd(s_1,N) = 1$ and set $t_1 := N /r_1$. Then for each $1 \leq k_1 \leq N$ holds
$$\left| \frac{p_1}{q_1} - \frac{k_1}{N} \right| \left| \frac{p_1t_1-k_1s_1}{Ns_1} \right| \geq \frac{1}{Ns_1},$$
which is the second claim of Theorem~\ref{main:thm} (i). This observation also explains entirely, what happens in Example~3.1 of \cite{FGW21}, where $\xi_1=\xi_2= \tfrac{1}{2}$.\\[12pt]
If $\xi_1,\xi_2 \notin \mathbb{Q}$, we again fix $N \in \mathbb{N}$ and only place weight at $x_1$ and $x_2$. This implies that the approximation error of the measure at $x_2$ is automatically equal to zero. Next, it is clear, that there exists exactly one $p(N) \in \mathbb{N}$ with $|\xi_1-p(N)/N| \leq 1/2N$. Multiplying the equation by $N$ we obtain $|N \xi- p(N)| \leq 1/2$. Hence, for arbitrary $N$ the best possible error term is $\norm{\left\{N \xi_1\right\}}$, where $\norm{x} : = \min(|x|,1-|x|)$. This means that the lower bound is governed by the distance of the Kronecker sequence $\left\{N \xi_1\right\}$ from its closest integer. As the Kronecker sequence is a uniformly distributed sequence, see e.g. \cite{Nie92}, Theorem~3.3, we thus have for every $c > 2$ that 
$$D_N^*(\mu,\nu_N) \geq \frac{1}{c N}$$
holds for infinitely many $N \in \mathbb{N}$. On the other hand, uniform distribution of the Kronecker sequence also implies
$$D_N^*(\mu,\nu_N) \leq \frac{1}{c N}$$
for infinitely $N \in \mathbb{N}$ because the approximation error at $x_1$ can get arbitrarily small.
\end{example}
Example~\ref{ex:n=2} also shows that the bounds in Theorem~\ref{main:thm} are sharp. In fact corresponding examples which proof the sharpness of our bounds can be constructed in the same manner for arbitrary $n \in \mathbb{N}$ (and also for infinitely supported discrete measures).\\[12pt]
Finally, we compare Theorem~\ref{main:thm} to classical results from Diophantine approximation: If the measure is supported on finitely many points $x_1,\ldots,x_n$ and the weights $\xi_i$ are linearly independent (over $\mathbb{Q}$), irrational, algebraic numbers, then Schmidt's subspace Theorem, see \cite{Sch72a}, can be applied and yields an error term of at least $N^{-(1+1/n+\epsilon)}$ for all $N \geq N_0$. If $N$ is big enough, this implies a lower bound for the star-discrepancy of order $N^{-(1+1/n+\epsilon)}$ which is worse than what we obtain. On the contrary, the simultaneous Dirchlet Theorem implies that it is possible to find infinitely $N$ such that each individual point mass is approximated of order $\leq \frac{1}{N^{1+1/n}}$. This is worse than the result $1/N$ from Theorem~\ref{thm:01} and thus does not impose an obstacle for what follows.

\section{Proof of the Main Result}
This section is dedicated to the proof of Theorem~\ref{main:thm}. 
\paragraph{Finite discrete measures.} At first, we consider part (i) of Theorem~\ref{main:thm}, i.e. a finite discrete measure with rational weights. We start with a straightforward remark.
\begin{remark} \label{rem:discrete} Let $\mu = p_1/q_1 \delta_x + (1-p_1/q_1) \tilde{\mu}$ be a discrete measure, where $\tilde{\mu}$ has no point mass at $x$. Assume that $q_1 | N$. Then the optimal approximative measure of $\mu$ assigns weight $p_1$ to $x$. 
\end{remark}
\begin{proof}[Proof of part (i)] If $\textrm{lcm}(q_1,\ldots,q_{n-1}) | N$, then $q_i|N$ for all $i=1,\ldots,n-1$. Thus we can assign multiplicities $k_i$ to the points $y_1=x_1,\ldots,y_n=x_n$ with $D_N^*(\mu,\nu_N) = 0$ by Remark~\ref{rem:discrete}. If $\textrm{lcm}(q_1,\ldots,q_{n-1}) \nmid N$, let $q_l$ be the smallest $q_i$ which does not divide $N$. By Remark~\ref{rem:discrete}, we can neglect all $i$ with $i<l$ and therefore, we may without loss of generality assume that $q_1 \nmid N$. We can (uniquely) write $q_i = r_i s_i$ with $r_i | N$ and $\gcd(s_i,N) = 1$ for each $i=1,\ldots,n$. Finally, set $t_1 = N /r_1$ and consider the interval $[0,x_1]$. Then for any measure associated to a finite point set
$$\mu([0,x_1]) - \nu_N([0,x_1]) = \left| \frac{p_1}{q_1} - \frac{k_1}{N} \right| = \left| \frac{p_1t_1-k_1s_1}{Ns_1} \right| \geq \frac{1}{Ns_1} \geq \frac{1}{Nq}$$
holds. The penultimate inequality holds because $\gcd(N,s_1) = 1$ yields $\gcd(t_1,s_1) = 1$ and $\gcd(p_1,q_1) = 1$ implies $\gcd(p_1,s_1) = 1$. 
\end{proof}
Next, we consider finite discrete measures with irrational weights $\xi_i$.
\begin{proof}[Proof of part (ii)] In comparison to Example~\ref{ex:n=2}, the situation can be treated as follows: at first, we take the best individual approximation of $\left\{N \xi_1\right\}$. Then we inductively define $p_l(N)$ by 
$$\left| \sum_{j=1}^l \xi_j - p_j(N)/N \right| < 1/2N,$$ 
As we assume $x_1 < x_2 < \ldots < x_n$, this algorithm guarantees to obtain the optimal discrepancy.\\[12pt]
Since $\sum_{i=1}^l \xi_i \notin \mathbb{Q}$ and the one-dimensional Kronecker sequence $\{N\sum_{i=1}^l \xi_i\}$ is uniformly distributed in $[0,1]$ (see e.g. \cite{KN74}, Section 2.3), also $\norm{ \{N \sum_{i=1}^l \xi_i\}}$  is uniformly distributed in $[0,1/2]$. Hence, for any $c > 2$, there exist infinitely many $N$ and finite point sets $x_1,\ldots,x_N$ with $D_N(\mu,\nu_N) \geq \frac{1}{c N}$ as in the case for $n=2$ points.
\end{proof}

\paragraph{Infinite discrete measures.} Finally, we come to the infinite case and make use of what we have already proven for finitely supported measures.
\begin{proof}[Proof of part (iii)]
Assume that $\mu = \sum_{i=1}^\infty \xi_i \delta_{x_i}$ and let $N \in \mathbb{N}$. Choose $z_1$ as the supremum over all $x \in [0,1]$ with $\mu([0,x]) \leq 1/2N$. At first let us assume that $\mu(z_1) = 0$. If we put weight at $z_1$ (or any point smaller than $z_1$), then $D_N^*(\mu,\nu_N) \geq 1/2N$ because the minimum weight we can choose is $1/N$. If we do not put weight at $z_1$, then the interval $[0,x_1]$ does not contain a point and we also have $D_N^*(\mu,\nu_N) \geq 1/2N$. Therefore, we obtain the desired lower bound in this case (which is not surprising because it very much resembles the continuous case). So the remaining case to solve is $\mu(z_1) > 0$. Again putting weight left of $z_1$ would result in $D_N^*(\mu,\nu_N) \geq 1/2N$. Let $\tilde{c} > 0$ be arbitrary. As long as $0 < \mu(z_1) < 1/\tilde{c}N$ and we put weight at $z_1$ we have $D_N^*(\mu,\nu_N) \geq (1/2 - 1/\tilde{c})N$. We now fix $z^* = z_1$ and let $N$ be big enough such that $\mu(z^*) > 1/2N$.\\[12pt] 
Next we define $z_i$ for $i=1,\ldots N/2$ to be the supremum over all $x \in [0,1]$ with $\mu([0,x]) \geq (2i-1)/N$. Since $\mu(z^*) > 1/2N$, the point $z^*$ is contained in the set of the $z_i$. Let us now consider the measure $\mu([0,z^*]) = \xi$. Since we have to put weight on $z^*$ (because otherwise $D^*(\mu,\nu_N) > 1/2N$ and we would be done), we cannot approximate it any better than $\{N \xi\}$ and we can apply the result either for the rational case if $\xi$ is rational or for the irrational case else. In any case, the claim follows. 
\end{proof} 

\bibliographystyle{alpha}
\bibdata{references}
\bibliography{references}

\textsc{Ruhr West University of Applied Sciences, Duisburger Str. 100, D-45479 M\"ulheim an der Ruhr,} \texttt{christian.weiss@hs-ruhrwest.de}\\[12pt]
\textsc{Max-Planck-Institut f\"ur Mathematik, Vivatsgasse 7, D-53111 Bonn,} \texttt{cweiss@mpim-bonn.mpg.de}

\end{document}